\documentclass[10 pt]{amsart}
\usepackage{amsfonts}
\usepackage{amssymb}
\usepackage{graphicx}

\allowdisplaybreaks
\let\oldmarginpar\marginpar
\renewcommand\marginpar[1]{\-\oldmarginpar[\raggedleft\footnotesize #1]%
{\raggedright\footnotesize #1}}

\def\({\left(}       \def\){\right)}
    
\def\U{\mathrm{U}}

\newtheorem{prop}{\sc Proposition}

\newtheorem{lem}{\sc Lemma}
\newtheorem{thm}{\sc Theorem}
\newtheorem{cor}{\sc Corollary}

\newenvironment{pf}{\noindent{\textit{Proof. }}}{$\Box$ }


\begin{document}
\title[Continuous deformations of harmonic maps]
{Continuous deformations of harmonic maps and their unitons}

\author[A. Aleman]{Alexandru Aleman}
\address{Department of Mathematics, University of Lund
Box 118, 22100 Lund, Sweden}\email{alexandru.aleman@math.lu.se}

\author[M. J. Mart\'in]{Mar\'ia J. Mart\'in}
\address{Department of Physics and Mathematics, University of Eastern Finland, P.O. Box 111, FI-80101 Joensuu, Finland.} \email{maria.martin@uef.fi}

\author[A. M. Persson]{Anna-Maria Persson}
\address{Department of Mathematics, University of Lund
	Box 118, 22100 Lund, Sweden}\email{Anna-Maria.Persson@math.lu.se}

\author[M. Svensson]{Martin Svensson}
\address{Department of Mathematics and Computer Science, University of Southern Denmark, 5230 Odense~M, Denmark}\email{svensson@imada.sdu.dk}

\thanks{The first and second authors are partially supported by MINECO/FEDER-EU research project MTM2015-65792-P, Spain. The second author also thankfully acknowledges partial support from Academy of Finland grant $268009$ and Spanish Thematic ResearchNetwork MTM\-2015-69323-REDT, MINECO}

%%%%%%%%%%%%%%%%%% title, etc.
\begin{abstract} We consider harmonic maps on simply connected Riemann surfaces into the group $\U(n)$ of unitary matrices of order $n$. It is known that a harmonic map with an associated algebraic extended solution can be deformed into a new harmonic map that has an $S^1$-invariant associated extended solution. We study this deformation in detail and show that the corresponding unitons are smooth functions of the deformation parameter and real analytic along any line through the origin.
\end{abstract}
\maketitle
%%%%%%%%%%%%%%%%%%%%%%%%% INTRODUCTION. %%%%%%%%%%%%%%%%%%%%%%%%%
\section*{Introduction}
%%%%%%%%%%%%%%%%%%%%%%%%%%%%%%%%%%%%%%%%%%%%%%%%%%%%%%%%%%%%%%%%%
\label{intro}
A smooth map $\varphi$ between two Riemannian manifolds $(M,g)$ and $(N,h)$ is said to be \emph{harmonic} if it is a critical point of the energy functionals
\[
E(\varphi, D)=\frac{1}{2}\int_D|\mathrm{d}\varphi|^2\omega_g\,,
\]
where $D$ is relatively compact in $M$, $\omega_g$ is the volume measure, and $|\mathrm{d}\varphi|^2 $ is the Hilbert-Schmidt  norm of the differential  of $\varphi$.
Using coordinates it is easy to see that
 $E(\varphi, D)$ is the natural  generalization of the classical integral of Dirichlet in $\mathbb{R}^d$
\[
D[u]=\int_\Omega |\nabla u|^2\, \mathrm{d} V\qquad (\Omega\subset \mathbb{R}^d)\,.
\]

In this paper we consider harmonic maps from a simply connected Riemann surface $M$ into the group $\mathrm{U}(n)$ of unitary matrices of order $n$. Such maps are interesting for many reasons. For example (see \cite{SW}), when $M$ coincides with the Riemann sphere $S^2$,  they are equivalent to harmonic maps of finite energy from the plane,  they  provide a nonlinear $\sigma$-model for particle physics, and  they give minimal branch immersions of $S^2$. As any compact Lie group has a totally geodesic embedding into $\mathrm{U}(n)$ for some $n$, these maps include harmonic maps into compact Lie groups and, in particular, harmonic maps into symmetric spaces. There is a vast literature on the subject (for example, we refer the reader to \cite{SW} and the references therein)  and some of the existent approaches combine ideas from differential geometry and operator theory.

The aim of this paper is to use Blaschke-Potapov factorizations of bounded analytic matrix-valued functions in the unit disc in order to prove  a deformation result for an important class of such maps which will be briefly described below.

Given a smooth map $\varphi:M\to \mathrm{U}(n)$  one  can use  a standard variational argument to derive an equivalent differential condition for harmonicity. More precisely (see \cite{U} or \cite[Sec. 3.1]{SW}), consider the matrix-valued 1-form
\[
\frac{1}{2}\varphi^{-1}\mathrm{d}\varphi=A^\varphi_z \mathrm{d} z+A^\varphi_{\bar{z}}\mathrm{d}\bar{z}\,,
\]
where $z$ is a  local coordinate on $M$. Then it turns out that  $\varphi$  is harmonic if and only if
\begin{equation}\label{harmcond}
(A^\varphi_z)_{\bar{z}}+(A^\varphi_{\bar{z}})_z=0\,.
\end{equation}

In her study of harmonic maps \cite{U}, Uhlenbeck introduced the notion of an \emph{extended solution}, which is a map $\Phi:S^1\times M\to \mathrm{U}(n)$ satisfying $\Phi(1,\cdot)=I$ and such that, for every local coordinate $z$ on $M$, there are $\mathfrak{gl}(\mathbb{C}^n)$-valued maps $A_z$ and $A_{\bar{z}}$ for which
\begin{equation}\label{extsol}
\Phi(\lambda,\cdot)^{-1}\mathrm{d}\Phi(\lambda,\cdot)=(1-\lambda^{-1})A_z\mathrm{d} z+(1-\lambda)A_{\bar{z}}\mathrm{d}\bar{z}\,.
\end{equation}
In this case $\varphi=\Phi(-1,\cdot)$ is a harmonic map with $A^\varphi_z=A_z$ and $A^\varphi_{\bar{z}}=A_{\bar{z}}$.

For a given harmonic map $\varphi:M\to \mathrm{U}(n)$, an extended solution with the property that
\[
\Phi^{-1}(\lambda,\cdot)\mathrm{d}\Phi(\lambda,\cdot)=(1-\lambda^{-1})A^\varphi_z\mathrm{d} z+(1-\lambda)A^\varphi_{\bar{z}}\mathrm{d}\bar{z}
\]
is said to be \emph{associated} to $\varphi$ and we have
\[
\Phi(-1,\cdot)=u\varphi
\]
for some constant $u\in \mathrm{U}(n)$.

As $M$ is assumed simply connected, the existence of extended solutions is granted by (actually, equivalent to) the harmonicity condition  (\ref{harmcond}). Moreover,  as a function of $z\in M$,  $\Phi(\lambda,\cdot)$ is as smooth as $\varphi$. The extended solution is called {\it algebraic} if  it is a trigonometric polynomial, i.e., there exist $r,s \in \mathbb{N}$ such that
\[
\Phi(\lambda,z)=\sum_{k=-r}^s T_k(z)\lambda^k\,.
\]

One of the  remarkable results in \cite{U} is that, for the Riemann sphere $M=S^2$, every harmonic map has an algebraic extended solution (in fact, any extended solution on a \emph{compact} Riemann surface is algebraic up to left multiplication by a constant loop, see \cite{OV}). Now standard factorization theory of matrix-valued functions on  $S^1$ (see \cite{Peller}) shows that such functions are essentially polynomial  Blaschke-Potapov products depending on $z\in M$.  More precisely,  there exist  subbundles $\alpha_1,\ldots,\alpha_m$ of $M\times\mathbb{C}^n$  such that
\begin{equation}\label{bp}
\Phi(\lambda,z)=\lambda^{-r}\prod_{j=1}^{m}(\pi_{\alpha_j(z)}+\lambda\pi_{\alpha_j(z)}^\perp),
\end{equation}
where $\pi_\alpha$ denotes the orthogonal projection onto the subspace $\alpha$ of $\mathbb{C}^n$ and $\pi_\alpha^\perp=\pi_{\alpha^\perp}$. It is well known that for  $n>1$,  the factors in the product above are not necessarily unique. This  is one of the major technical  difficulties in dealing with these objects. However, it is shown in \cite{U} that the factors can be chosen such that for $k \le m$ the partial products
$$
\Phi_k(\lambda,z)=\lambda^{-k}\prod_{j=1}^{k}(\pi_{\alpha_j(z)}+\lambda\pi_{\alpha_j(z)}^\perp)
$$
are  extended solutions as well. The Blaschke-Potapov factors involved in such a factorization are called {\it unitons} of   $\Phi$ and (\ref{bp}) is its {\it uniton factorization}.
For $\lambda=-1$ this yields a factorization of the original harmonic map $\varphi$ also called a uniton factorization,  and the corresponding factors are called unitons of $\varphi$.  A harmonic map which admits such a factorization or, equivalently, has an associated algebraic extended solution, is said to have {\it finite uniton number}.

An extended solution $\Phi$ is said to be {\it $S^1$-invariant}  (see \cite[Sec. 10]{U} or \cite[Secs. 2.3 and 3.3]{SW}) if it satisfies
\[
\Phi(\lambda_1\lambda_2,z)=\Phi(\lambda_1,z)\Phi(\lambda_2,z)\qquad(z\in M,~\lambda_1,\lambda_2\in S^1)\,.
\]

If $\Phi$ is an algebraic $S^1$-invariant extended solution, there exist integers $k_1,\dots,k_n$ and a loop $\gamma:S^1\to \mathrm{U}(n)$ of the form
\[
\gamma(\lambda)={\rm diag}(\lambda^{k_1},\lambda^{k_2},\dots,\lambda^{k_n})
\]
such that $\Phi$ takes values in the orbit of $\gamma$ under the adjoint action of $\mathrm{U}(n)$. Equivalently, $\Phi$ has a uniton factorization (\ref{bp}) where the unitons satisfy
\[
\alpha_j\subset\alpha_{j+1},\quad \partial_z\alpha_j\subset\alpha_{j+1}\qquad(j=1,\dots,k-1)
\]
and
\[
\partial_{\bar{z}}\alpha_j\subset\alpha_j\qquad(j=1,\dots,k)\,.
\]

Using an approach based on Morse theory, Burstall and Guest showed in \cite{BG} that an algebraic extended solution on $S^2$ can be deformed into an $S^1$-invariant extended solution. From this follows that a harmonic map from $S^2$ into a compact Lie group can be deformed through a family of harmonic maps to one that arises from a twistor construction in the sense of Burstall and Rawnsley \cite{BR}. As pointed out in the abstract, we will continue the study of this deformation and the dependence of the corresponding maps and their unitons on the deformation parameter. We will use an approach based on functional analysis instead of Morse theory.

To describe the deformation, assume that $\Phi$ is an algebraic extended solution. We write it in the form (\ref{bp}) and denote by $b$  the corresponding Blaschke-Potapov product, i.e.,
\[
b(\lambda,z)=\prod_{j=1}^{m}(\pi_{\alpha_j(z)}+\lambda\pi_{\alpha_j(z)}^\perp).
\]
Then for fixed $\mu\in \mathbb{C}\setminus\{0\}$ and $z\in M$, $\lambda\mapsto b (\mu\lambda,z)$ is a matrix-valued polynomial whose determinant vanishes only at the origin. Therefore (see \cite{Peller}) it can be factored in the form $b^\mu G^\mu$, where $G^\mu$ is analytic and invertible in a neighborhood of the unit disc and
\begin{equation}\label{bdeformation}
b^\mu(\lambda,z)=\prod_{j=1}^{m}(\pi_{\alpha_j(\mu,z)}+\lambda\pi_{\alpha_j(\mu,z)}^\perp)\,.
\end{equation}
The deformation of $\Phi$ is then given by the family $\{\Phi^\mu\}$ where
\begin{equation}\label{deformation}
\Phi^\mu(\lambda,z)=\lambda^{-r}b^\mu(\lambda,z)\,.
\end{equation}
Clearly, $\Phi^1=\Phi$.

We are going to prove the following theorem generalising results by Burstall and Guest \cite{BG}.

\begin{thm}\label{main} Let $\varphi:M\to \mathrm{U}(n)$ be a harmonic map with algebraic extended solution $\Phi$ and let $\Phi^\mu$ be its deformation defined by (\ref{deformation}).
\begin{itemize}
\item[(i)] For each $\lambda\in S^1$ the function $(\mu,z)\to \Phi^\mu(\lambda,z)$ extends to $\mathbb{C}\times M$,  is $C^\infty$ in  the variable $\mu\in \mathbb{C}$, and its restriction to any line through the origin is real-analytic in $\mu$.  Moreover, for each $\mu\in \mathbb{C}$, $\Phi^\mu$ is an algebraic extended  solution.
\item[(ii)] $\Phi^0(-1,\cdot)$  is $S^1$-invariant.
\item[(iii)] Given any uniton factorization of  the form (\ref{bp}) of $\Phi$, there exists for each $\mu\in \mathbb{C}$ a uniton factorization of $\Phi^\mu$,
\[
\Phi^\mu(\lambda,z)=\lambda^{-r}\prod_{j=1}^{m}(\pi_{\alpha_j(\mu,z)}+\lambda\pi_{\alpha_j(\mu,z)}^\perp)\,,
\]
such that for $1\le j\le m$, the uniton $b_j^\mu(\lambda,z)=\pi_{\alpha_j(\mu,z)}+\lambda\pi_{\alpha_j(\mu,z)}^\perp$ is a $C^\infty$-function of $\mu$ which is real-analytic on any line through the origin and satisfies $b^1_j(\lambda,z)=b_j(\lambda,z)$.
\end{itemize}
\end{thm}
\section{Preliminaries}

Denote by $\mathcal{H}_+$ the usual Hardy space of $\mathbb{C}^n$-valued functions identified with the closed subspace of  $L^2(S^1,\mathbb{C}^n)$ consisting of Fourier series whose negative coefficients vanish.

We shall make use of a method originating from the theory of integrable systems, called the Grassmanian model \cite{Se},  which associates to an extended solution $\Phi$ the family of closed subspaces $W(z),~z\in M,$   of  $L^2(S^1,\mathbb{C}^n)$, defined by
\[
W(z)=\Phi(\cdot,z)\mathcal{H}_+\,.
\]
Let us denote by $\partial_z$ and $\partial_{\overline{z}}$ the derivatives with respect to  $z$ and $\overline{z}$ on $M$, respectively, and by $T$ the forward shift on $L^2(S^1,\mathbb{C}^n)$:
\[
(Tf)(\lambda)=\lambda f(\lambda)\qquad(\lambda \in S^1)\,.
\]

If $f:S^1\times M\to\mathbb{C}^n$ is differentiable in the second variable and satisfies $f(\cdot,z)\in W(z),~z\in M$, it follows from (\ref{extsol}) that $T\partial_zf(\cdot,z)\in W(z)$  and $\partial_{\overline{z}}f(\cdot,z)\in W(z)$, i.e., in terms of  differentiable sections we have
\begin{equation}\label{propW}
T\partial_z W(z)\subset W(z),\quad \partial_{\overline{z}}W(z)\subset W(z)\,.
\end{equation}
This technique is relatively common in the theory of integrable systems and has the advantage that it yields coordinate-free equations, since the derivatives of the coordinates involved in the chain rule  are absorbed in the corresponding subspaces. The converse of this statement is important for our purposes. A more general version of the following result can be found in \cite{Se}.

\begin{prop}\label{W-char}
Let $\Phi:S^1\times M\to \mathrm{U}(n)$ be smooth with $\Phi(1,z)=I,~z\in M$. If $W(z)=\Phi(\cdot,z)\mathcal{H}_+$ satisfies (\ref{propW}), then $\Phi$ is an extended solution, i.e., it satisfies the equation (\ref{extsol}) and $\Phi(-1,\cdot)$ is harmonic.
\end{prop}

There is another simple observation regarding smoothness (in $z$) related to such families of shift-invariant subspaces. Suppose that $N$ is a smooth manifold and that we have a function from $N$ into the set of polynomial Blaschke-Potapov products, $N \ni \nu\to b^\nu,$
and let $W_\nu=b^\nu\mathcal{H}_+$.  Then obviously, their orthogonal complements $W_\nu^\perp$ in $\mathcal{H}_+$ have finite dimension.

\begin{lem}\label{smoothness} $W_\nu^\perp$ form a $C^k$ ($C^\omega$) vector bundle over $N$ if and only if for every $\lambda\in S^1$, the map $\nu\to b^\nu(\lambda)$ is $C^k$ ($C^\omega$).\end{lem}
\begin{pf} The connection between the two objects is given by the reproducing kernels in $W_\nu^\perp$. For example, if $\zeta$ belongs to the unit disc and $\{e_1,\ldots, e_n\}$ is the canonical basis in $\mathbb{C}^n$, the functions
\begin{equation}\label{reprokernel}
k_\zeta^i(\lambda,\nu)=(1-\overline{\zeta}\lambda)^{-1}(I-b^\nu(\lambda)(b^\nu)^*(\zeta))e_i
\end{equation}
satisfy
\[
\langle f,k_\zeta^i(\cdot,\nu)\rangle_{\mathcal{H}_+}=\langle f(\zeta),e_i\rangle_{\mathbb{C}^n}\qquad (f\in W_\nu^\perp)\,.
\]
In particular, if $\{p_0(\cdot,\nu),\ldots,p_l(\cdot,\nu)\}$ is an orthonormal basis in $W_\nu^\perp$ then
\begin{equation}\label{kernelid}
k_\zeta^i(\lambda,\nu)=\sum_{k=0}^l \langle e_i, p_k(\zeta,\nu)\rangle_{\mathbb{C}^n}\, p_k(\lambda,\nu)\,.
\end{equation}	
Obviously, the span of these reproducing kernels is $W_\nu^\perp$. Thus if  the maps $\nu\to b^\nu(\lambda)$ are $C^k$ ($C^\omega$) we can construct bases in $W_\nu^\perp$ which are locally $C^k$ ($C^\omega$), and conversely, if $W_\nu^\perp$ from a $C^k$ ($C^\omega)$ vector bundle over $M$ we can find locally near every point an orthonormal basis whose elements are $C^k$ ($C^\omega$) functions of $\nu$ and by (\ref{kernelid}) and (\ref{reprokernel})  it follows that  $\nu\to b^\nu(\lambda)$ is $C^k$ ($C^\omega$).
\end{pf}

As an immediate application of the above results we obtain that the functions $\Phi^\mu,~\mu\in \mathbb{C}\setminus\{0\}$,  given in Theorem \ref{main} are smooth extended solutions.
\begin{cor}\label{def-extsol}  For $\mu \in\mathbb{C}\setminus  \{0\}$,  the function $\Phi^\mu$ defined by (\ref{deformation}) is an algebraic extended solution which is a smooth function of $z\in M$.
\end{cor}
\begin{pf} If $\Psi^\mu(\lambda,z)=\Phi(\mu\lambda,z)$, by the remarks in the Introduction,  we obviously have that
\[
W(z)^\mu=\Phi^\mu(\cdot,z)\mathcal{H}_+=\Psi^\mu(\cdot, z)\mathcal{H}_+.
\]
Then $W(z)^\mu$ satisfy (\ref{propW}) and by Proposition \ref{W-char} it is an extended solution which, by definition, is algebraic.
To verify the smoothness in $z\in M$, recall the notation $b^\mu(\lambda,z)=\lambda^r\Phi^\mu(\lambda,z)$ and note that the functions in $(b^\mu(\cdot,z)\mathcal{H}_+)^\perp$ are polynomials $p$ of degree at most $m-1$ whose coefficients, $\hat{p}_0,\ldots\hat{p}_{m-1}$, are solutions of the homogeneous linear system
\begin{equation}
\label{parseval}
\sum_{k=0}^{m-j-1}\langle \hat{p}_{k+j}, b^1_ke_i\rangle \mu^k=0 \qquad (0\le j\le m-1\,,\ 1\le i\le n),
\end{equation}
where $b^1_k$ are the Fourier coefficients of $b^1=T^r\Phi$. Then $(b^\mu(\cdot,z)\mathcal{H}_+)^\perp$ is a smooth bundle over $M$ and the result follows by Lemma \ref{smoothness}.
\end{pf}

The above argument provides a smooth dependence of $\Phi^\mu$ of the variable $\mu\in \mathbb{C}\setminus\{0\}$ but not at the origin, since the linear system (\ref{parseval}) in the proof of Corollary~\ref{def-extsol} degenerates. In order to prove our main theorem we need to overcome this difficulty.

Our second preliminary observation concerns the representation of  the  Blaschke-Potapov products involved in our considerations. These are matrix-valued polynomials which are unitary on $S^1$ and whose determinants vanish only at the origin. As pointed out in the Introduction, such a function $b$ can be written in the form
\begin{equation}
\label{bp1}
b(\lambda)=\prod_{j=1}^{m}(\pi_{\alpha_j}+\lambda\pi_{\alpha_j}^\perp)\,,
\end{equation}
where $\alpha_1,\ldots,\alpha_m$ are subspaces of $\mathbb{C}^n$ and $\pi_{\alpha_j}$ are the corresponding orthogonal projections. The factorization is not unique and in the case of extended solutions, not every factorization yields unitons. This situation is analyzed in  \cite{SW}, where factorizations are obtained with help of {\it filtrations} as follows.

Given a shift-invariant subspace $W$ of $\mathcal{H}_+$, a filtration associated to it is a  finite sequence of shift invariant subspaces $W_k,~0\le k\le m,$ such that $W_m=W,~W_0=\mathcal{H}_+$, and
$$
TW_{k-1}\subset W_k\subset W_{k-1}\qquad(1\le k\le m).
$$
Clearly, every factorization of the form (\ref{bp1}) of a polynomial Blanchke-Potapov product $b$ yields the filtration given by
\begin{equation}\label{filtration}
W_k=\prod_{j=1}^{k}(\pi_{\alpha_j}+\lambda\pi_{\alpha_j}^\perp)\mathcal{H}_+\qquad(1\le k\le m)\,.
\end{equation}
Conversely, it is not difficult to verify that every filtration of $b\mathcal{H}_+$ gives a factorization of the form (\ref{bp1}). A useful example, referred to in \cite{SW} as the \emph{Segal filtration}, is given by
\[
W_k=W+T^{k}\mathcal{H}_+\qquad (0\le k\le m)\,,
\]
where $m$ is the degree of $b$. It yields  the standard factorization as described in \cite{Peller} (see Lemma~5.1, p. 76). Moreover, by definition we see that the intermediate spaces satisfy (\ref{propW}), so that for extended solutions, the procedure yields a uniton factorization. Also note that for $k\le m$, $W_k^\perp$ consists of the polynomials in $W^\perp$ whose degree does not exceed $k-1$. For our purposes, we shall use this filtration in order to exhibit a special type of basis in the orthogonal complement of our spaces $W(z),~z\in M$.
\begin{prop}\label{trianglebasis}
Let $\Phi$ be a polynomial extended solution and set $W=\Phi\mathcal{H}_+$. Then for every $z_0\in M$ there exists a neighborhood $V$ of $z_0$ such that for all $z\in V$ there is a basis
\[
 \mathcal{B}_z=\{p_0(\cdot,z),\ldots p_l(\cdot,z) \}
\]
of $W(z)^\perp$  consisting of polynomials that are smooth functions of $z\in M$,  with the properties:
\begin{itemize}
\item[(i)]The degree of $p_k(\cdot,z)$ is non-decreasing in $k$.
\item[(ii)] For  $j\ge 0$, the set of polynomials of degree at most $j$ in $\mathcal{B}$  spans the space of  polynomials of degree at most $j$ in  $W(z)^\perp$.
\end{itemize}
\end{prop}

\begin{pf}
We can construct $\mathcal{B}$ inductively as follows. Let $W_k,~0\le k\le m$, be the Segal filtration for $W$. Recall that for $k\ge 1$, $W_k^\perp$  consists of polynomials of degree at most $k-1$. It is also easy to verify that $W_k^\perp$ is a smooth subbundle of  the trivial bundle given by the space of polynomials whose degree does not exceed $k-1$. Start with $k=1$ and for $z$ near $z_0$  construct  linearly independent sets
$\mathcal{B}_z^0$,  smooth in $z$,    satisfying (i) and (ii) with $j=0$. Now assume that for $z$ near $z_0$  we have constructed  linearly independent sets  $\mathcal{B}_z^s$,  smooth in $z$, and satisfying (i) and (ii) with $j\le s$. Then
\[
{\rm span}\, \mathcal{B}_z^s =  (W(z)+T^{s+1}\mathcal{H}_+)^\perp\,.
\]

Assume that $s+2<m$   and consider the previous step in the Segal filtration, i.e. $W+T^{s+2}\mathcal{H}_+$. As pointed out above, $(W+T^{s+2}\mathcal{H}_+)^\perp$  consists of the polynomials in $W^\perp$ whose degree does not exceed $s+1$ and
if the  factorization of $\Phi$ corresponding to this filtration is
\[
 \Phi(\lambda,z)=\prod_{j=1}^{m}(\pi_{\alpha_j(z)}+\lambda\pi_{\alpha_j(z)}^\perp)\,,
\]
then
\[
W+T^{s+2}\mathcal{H}_+=b_s\mathcal{H}_+\,,
\]
with
\[
b_s(\lambda,z)= \prod_{j=1}^{s+2}(\pi_{\alpha_j(z)}+\lambda\pi_{\alpha_j(z)}^\perp)\,.
\]
As pointed out above, $(W+T^{s+2}\mathcal{H}_+)^\perp$ is a smooth bundle on $M$, hence by Lemma \ref{smoothness}, $b(\cdot,z)$ is smooth in $z$ and thus the reproducing kernels given in (\ref{reprokernel}),
\[
k_\zeta^i(\lambda,z)=(1-\overline{\zeta}\lambda)^{-1}(I-b_s(\lambda,z)b_s^*(\zeta,z))e_i\,,
\]
with $ \{e_i, 1\le i\le n\}$  the canonical basis in $\mathbb{C}^n$, and, say,  $0<|\zeta|<1$, are smooth as well and span this space.	
Then we can obviously find a possibly smaller neighborhood $\tilde{V}$  of $z_0$ and a finite set $F\subset \{1,\ldots,n\}\times \{\zeta:~|\zeta|<1\}$ such that in that neighborhood
\[
\mathcal{B}_z^{s+1}=\mathcal{B}_z^s\cup\{k_\zeta^i(\cdot,z):~(i,\zeta)\in F\}\qquad(z\in \tilde{V})
\]
satisfies (i) and (ii) for $j\le s+1$,  which completes the proof.
\end{pf}

\section{Proof of the main result}

Let $\varphi:M\to \mathrm{U}(n)$ be harmonic with algebraic extended solution $\Phi$. We shall assume throughout that $\Phi\mathcal{H}_+\subset\mathcal{H}_+$ since factors of the form $\lambda^{-r}$ play no role in our considerations. Thus the deformation $\Phi^\mu$ has the form (\ref{bdeformation}) and $W^\mu=\Phi^\mu\mathcal{H}_+\subset\mathcal{H}_+$.

(i) Let $\mu\in \mathbb{C}\setminus\{0\}$. The  simple argument based on Parseval's formula   which  gives (\ref{parseval}) (proof of Corollary~\ref{def-extsol}) shows that a polynomial $p$ belongs to $W(z)^\perp$ if and only if the polynomial
\[
p_\mu(\lambda)=p\left(\frac{\lambda}{\mu}\right)
\]
belongs to $(W(z)^\mu)^\perp$. Now let $z_0\in M$ be arbitrary and consider smooth  bases
\[
\mathcal{B}_z=\{p_0(\cdot,z),\ldots p_l(\cdot,z) \}
\]
of $W(z)^\perp$, defined locally near $z_0$, with the  properties (i) and (ii) in Proposition~\ref{trianglebasis}.
Set
\[
q_k^\mu(\lambda,z)= \mu^{d_k}(p_k)_\mu(\lambda,z)=\mu^{d_k} p_k\left(\frac{\lambda}{\mu},z\right)\,,
\]
where $d_k$ is the degree of $p_k$, and note that $\{q_k^\mu(\cdot,z):~0\le k\le l\}$ is a basis in $(W(z)^\mu)^\perp$. The entries of the matrix
\[
\mathcal{A}(\mu,z)=(\langle q_i^\mu(\cdot,z),q_j^\mu(\cdot,z)\rangle_{\mathcal{H}_+})_{0\le i,j\le l}
\]
are polynomials in $\mu$ and $\overline{\mu}$. Thus, $\mathcal{A}(\cdot,z)$ extends to a $C^\infty$-function on $\mathbb{C}$ which is real-analytic on any line through the origin since on such lines the entries of $\mathcal{A}$ become analytic polynomials in one variable. Moreover, $\det\mathcal{A}(\mu,z)\ne 0$ for $\mu\in \mathbb{C}\setminus \{0\}$ and $z$ near $z_0$.  We claim that
\[\det\mathcal{A}(0,z)\ne 0\,,
\]
for $z$ near $z_0$.

Indeed, if we assume the contrary, it follows easily by the definition of $q_k^\mu(\cdot,z)$ that there exist polynomials in $\mathcal{B}_z$ with the same degree such that their dominant terms are linearly dependent. But this means that the same linear combination of these polynomials has a strictly smaller degree, which contradicts Proposition \ref{trianglebasis} (ii).

Now write each element $p\in (W(z)^\mu)^\perp$  as a linear combination
\[
p=\sum_{j=0}^lc_j(p,\mu,z)q_j^\mu(\cdot,z)
\]
and take scalar products with $q_k^\mu(\cdot,z)$ in $\mathcal{H}_+$ on both sides in order to conclude that the coefficients are given by
\[
c_j(p,\mu,z)=\mathcal{A}^{-1}(\mu,z)(\langle p,q_j^\mu(\cdot,z)\rangle_{\mathcal{H}_+})\,.
\]
In the case of the reproducing kernels given by (\ref{reprokernel}), the vector on the right hand side is $(\langle e_i,q_j^\mu(\zeta,z)\rangle$, hence by Cramer's rule and the above argument it follows easily that $c_j(k_\zeta^i, \mu,z)$ are  $C^\infty$-functions of $\mu\in\mathbb{C}$, which are real-analytic on any line through the origin, hence so is each $k_\zeta^i$. Since our spaces consist of polynomials, we can choose $\zeta=1$, and from the fact that $\Phi^\mu(1,z)=I$ it follows that $\Phi^\mu(\lambda, z)e_i$ extends to a $C^\infty$-function on $\mathbb{C}$ which is real-analytic on any line through the origin.  By Corollary \ref{def-extsol}  we know that $\Phi^\mu$ is an algebraic solution, and using also the above argument we can conclude that $(\mu,z)\to\Phi^\mu(\lambda,z)$ has continuous derivatives in $\mathbb{C}\times M$. Then $\Phi^0$ satisfies (\ref{extsol}) with matrices independent of $\lambda$; in other words, it is an extended solution. Finally, $\Phi^0$ is a polynomial in $\lambda$ of the same degree as $\Phi^\mu,~\mu\in \mathbb{C}\setminus\{0\}$.

(ii)  Recall that $\Phi(\mu\lambda,z)=\Phi^\mu(\lambda,z)G^\mu(\lambda,z)$ with $G^\mu(\cdot,z)$ analytic and invertible in a neighborhood of the closed unit disc.  If $\mu_1,\mu_2\in S^1$ then
\[
\Phi^{\mu\mu_1}(\lambda,z)G^{\mu\mu_1}(\lambda,z)=\Phi(\mu\mu_1\lambda,z)=\Phi^{\mu}(\mu_1\lambda,z)G^\mu(\mu_1\lambda,z)
\]
and a standard argument gives
\[
\Phi^{\mu\mu_1}(\lambda,z)=\Phi^{\mu}(\mu_1\lambda,z)\Phi^\mu(\mu_1,z)^{-1}\,.
\]
Let $\lambda=\mu_2$ and let $\mu\to 0$ to obtain
\[
\Phi^0(\mu_2,z) \Phi^0(\mu_1,z)=\Phi^0(\mu_1\mu_2,z)\,,
\]
i.e., $\Phi^0(-1,\cdot)$ is $S^1$-invariant.

(iii)  Consider a uniton factorization of the form (\ref{bp}) of $\Phi$ and let $W_k,~0\le k\le m$, be the filtration (\ref{filtration}) induced by it.  Since this is a uniton factorization,  we have
\[
W_k=\Phi_k\mathcal{H}_+\,,
\]
with  extended solutions $\Phi_k$. It is easy to verify that if $\Psi_k^\mu(\lambda,z)=\Phi_k(\mu\lambda,z)$, then
\[
W_k^{\mu}=\Psi_k^\mu \mathcal{H}_+ \qquad (0\le k\le m)
\]
is a filtration with $W_0^\mu=W^\mu$. By the previous arguments, $$W_k^\mu=\Phi_k^\mu\mathcal{H}_+,$$
where $\mu\to \Phi_k^\mu$ extends to a  $C^\infty$-function on $\mathbb{C}$ which is analytic on each line through the origin. The unitons are given by
\[
b_k^\mu(\lambda,z)=\Phi_k^\mu(\lambda,z)\Phi_{k+1}^\mu(\lambda,z)^{-1}\,.
\]
For $\mu=1$ they obviously coincide with the original ones. Moreover, if $\lambda\in S^1$, $(\Phi_k^\mu)^{-1}=(\Phi_k^\mu)^*$, hence on $S^1$, $(\mu,z)\to b_k^\mu(\lambda,z)$ obviously has the smoothness required in the statement and since these functions are linear in $\lambda$ the properties hold everywhere. The proof is now complete.

\end{document}